\documentclass[amstex,12pt]{article} % CHANGED BY JPZ
\usepackage{graphicx} % CHANGED BY JPZ
\usepackage{amscd} % CHANGED BY JPZ
\usepackage{amssymb} % CHANGED BY JPZ
\usepackage{amsthm} % CHANGED BY JPZ
\usepackage{amsmath} % CHANGED BY JPZ
\newtheorem{definition}{Definition}[section]
                          
\newtheorem{prop}[definition]{Proposition}

\newtheorem{teo}[definition]{Theorem}

\newcommand{\R}{I\!\!R}
\newcommand{\N}{I\!\!N}

\newcommand{\BB}{{\cal B}}
\newcommand{\V}{{\cal V}}
\newcommand{\no}{\noindent}

\newcommand{\sfb}{$\Sigma-$flow-box} 
\newcommand{\xc}{$X\in {\frak X}^r(M)$ }

\newcommand{\xcr}{$X\in {\frak X}^r(M)$, $1 \le r \le \infty,$ } 
\newcommand{\xcrv}{$X\in {\frak X}^r(M)$, $1 \le r \le \infty,$ } 
\newcommand{\xcrp}{$X\in {\frak X}^r(M)$, $1 \le r \le \infty.$ } 
\newcommand{\vte}{virtual orthogonal edge}%%{virtual transversal edge}
\newcommand{\vtes}{virtual orthogonal edges}%%{virtual transversal edges}

\newcommand{\bkt}{$\BB_K-$type}
\newcommand{\tc}{T-closable}
\newcommand{\esp}{   }
\newcommand{\iet}{{\it iet}}
\newcommand{\ov}{\overline}

\begin{document}

\title{\bf On $C^r-$closing for flows on 2-manifolds. }

\author{Carlos Gutierrez \\
e-mail: gutp@impa.br \\
IMPA. \, \, 
Estrada Dona Castorina 110, 
J. Bot\^anico, \\
22460-320, 
Rio de Janeiro, R.J., 
Brazil.}
\date{}

\maketitle

\begin{abstract}
For some full measure subset $\BB$ of the set of 
\iet's (i.e. interval exchange
transformations) the following is satisfied:

Let $X$ be a $C^r$, $1\le r\le \infty$, vector field,
with finitely many singularities, on a compact 
orientable surface
$M$. Given a nontrivial recurrent point $p\in M$  of $X,$ 
the holonomy map around $p$ is semi-conjugate to 
an \iet \esp
$E :[0,1) \to [0,1).$ If $E\in \BB$ then
there exists a $C^r$ vector field $Y$,
arbitrarily close to $X$, in the $C^r-$topology, such that
$Y$ has a closed trajectory passing through $p$.
\end{abstract}

\section{Introduction}

The open problem 
``$C^r$-closing lemma'' is stated as follows:

``Let $M$ be a smooth compact manifold, $r\ge 2$ be
an integer, $f\in \mbox{Diff}^r(M)$ (resp. $X\in {\frak X}^r(M)$)
and $p$ be a nonwandering point of $f$ (resp. of $X$).
There exists $g\in \mbox{Diff}^r(M)$ (resp. 
$Y\in {\frak X}^r(M)$ ) arbitrarily close to $f$
(resp. to $X$) in the $C^r-$topology so that $p$ is
a periodic point of $g$ (resp. of $Y$)''.

C. Pugh proved the $C^1-$closing lemma  [Pg1]. 
There are few
previous results when $r\ge 2:$
Gutierrez [Gu1] showed  
similar results to this paper when the manifold is the torus $T^2$. 
There are  negative answers:
Gutierrez [Gu3] proved that
if the perturbation is localized in a
small neighborhood of the nontrivial recurrent point, then
$C^2-$closing is not always possible.
C. Carroll's [Car]  proved that, even for flows with
finitely many singularities, $C^2-$closing by 
a twist-perturbation (supported in a cylinder)
is not always possible. 
Concerning hamiltonian flows, M. Herman [Her] has 
remarkable counter-examples to the
$C^r-$closing lemma. 
Within the context of
geodesic laminations, S. Aranson and 
E. Zhuzhoma anounced in 1988 [A-Z] 
the $C^r-$closing lemma for a class of flows on surfaces;
however, their proofs have not been published yet.
For basic definition the reader may consult [K-H].

\section{Statement of the results}

Throughout this article, $M$ will be a smooth, 
orientable, compact,
two manifold and $\chi$ will be its Euler characteristic. 
We shall denote by ${\frak X}^r(M)$ the
space of vector field of class $C^r$,  
$1\le r \le \infty,$  with the $C^r-$topology. 
The trajectory of \xc \esp passing through
$p\in M$ will be denoted by $\gamma_p$
The domain of definition of a map $S$ will be denoted by
$\mbox{DOM}(S).$  
Smooth segments on $M$ will be denoted and referred as (open, half-open,
closed) intervals. 

A bijective map $E:[0,1) \to [0,1)$ is said to be
an \iet, i.e. an {\it Interval Exchange Transformation \/}
(with $m$ intervals)
if there exists a finite sequence
$0=a_1 < a_2 < \cdots < a_m < a_{m+1} = 1$ such that,
for all $i\in \{1,2,\cdots,m\}$ and for all $x\in [a_i, a_{i+1}),$ 
$E(x) = E(a_i) + x - a_i,$ and moreover, $E$ is discontinuous at exactly
$a_2, a_3, \cdots, a_m.$ This  $E$ will  be identified with 
the pair $(\lambda, \pi)\in \Delta_m\times {\frak S}_m$ 
made up of the positive
probability vector $\lambda = \{ |a_{i+1}- a_i| \}_{i=1}^m$ 
and the permutation
$\pi$ on the symbols $1, 2, \cdots, m,$ defined by
$\pi(i) = \#\{ j : E(a_j) \le E(a_i)\}.$ 
The space  of \iet's, with $m$ intervals, 
defined in $[0,1),$ 
will be identified with 
the {\it measurable space \/} $\Delta_m \times {\frak S}_m$ endowed  
with the product measure, where
$\Delta_m$ is the simplex of positive probability vectors of
$\R^m,$ with Lebesgue measure,
and ${\frak S}_m$ is the finite set of permutations
on $m$ symbols with counting measure. 
Let $E : [a,b) \to [a, b)$ be an \iet. We say that
$[s,t] \subset [a,b)$ is  
a {\it \vte \esp for $E,$ \/} if
$E$ restricted to $[s,t]$ is continuous and 
$s < E(s) < E^2(s) = t.$ Given $k\in \N,$ let $\BB_k$ be the set of \iet's \,  
$E : [a,b) \to [a, b)$ such that  for some sequence
$\, \, b_n \to a\, \, $ of points of $(a,b),$  
and for every $n\in \N,$ the \iet \esp \,
$E_n: [a,b_n) \to [a, b_n),$ induced by $E$, 
has at least $\chi+k+3,$ pairwise disjoint, \vtes. 
Denote $\BB = \cap_{k\ge 1}\BB_k.$
It will be seen that,
as a direct consequence of the work of 
W. A. Veech [Vee] and H. Masur [Mas],

\begin{teo} For all $m\ge 2,$ \,\, 
$\Delta_m \times {\frak S}_m \setminus \BB$
is a measure zero set.
\end{teo}

By transporting information along
flow boxes, Item (a2) below follows from the definition of $\BB_K.$ 

\begin{teo}([Gu2, Structure Theorem, Section 3])
Let $X\in {\frak X}^1(M).$
There are finitely many nontrivial recurrent trajectories
$\gamma_{p_1}, \gamma_{p_2}, \cdots, \gamma_{p_\ell}$ 
of $X$ such
that if $\gamma_p$ is any nontrivial recurrent trajectory of
$X,$ then $\ov{\gamma_p}= \ov{\gamma_{p_i}},$ 
for some $i=1,2,\cdots, \ell.$

Suppose that $X$ has exactly $K\in \N$ singularities
(K=0 is allowed).
Let $p\in M$ be a nontrivial recurrent point of $X.$ 
Take a half-open interval $[p,q)\subset M$  
transversal to $X,$
such that $p$ is a cluster point of 
$\gamma_p \cap (p,q),$
Denote by 
$P_X: [p,q) \to [p,q)$  the forward 
Poincar\'e map induced by $X$. 
If $[p,q)$ is small enough, it can be 
associated to $(p,[p,q)),$
an \iet \esp
$E=E_{(p,[p,q))}:[0,1) \to [0,1)$ 
and a 
continuous monotone surjective map
$h:[p,q) \to [0,1)$ such that
$\, h(p) = 0,\,$ 
$h$ restricted
to any given orbit of $P_X$ is injective
and, for all $x\in \mbox{DOM}(P_X),$ 
$\, \, E\circ h(x) = h \circ P_X(x);$ moreover,  
\begin{enumerate}
\item[{(a1)}]  there exists a subset 
$S\subset [0,1)$ of at most $\chi+K+2$ elements
such that if $A$ is a connected component of 
$[0,1)\setminus S,$
then $h^{-1}(A)$ is contained in $\mbox{DOM}(T);$
\item[{(a2)}] 
Let $\ov{p}\in \ov{\gamma_p}$ be a 
nontrivial recurrent point of $X$ and 
$(\ov{p}, [\ov{p}, \ov{q}))$
be a pair satisfying the same conditions 
as those of $[p, [p,q))$ above.
Then the property that the \iet \esp \, $E_{(\ov{p}, [\ov{p}, \ov{q}))}$
belongs to $\BB_K$ does not depend on $(\ov{p}, [\ov{p}, \ov{q})).$
\end{enumerate}
\end{teo}

Under conditions of theorem above and if $E\in \BB_K,$
any nontrivial recurrent point of $\ov{\gamma_p}$ is said to 
be of {\it \bkt.} 
{\it Our result is the combination of
Theorems 2.1 - 2.3.}

\begin{teo}
Let \xcr have $K\ge 0$ singularities.
Let $p\in M$ be a \bkt \esp nontrivial recurrent point of $X.$ 
Then
there exists  $Y\in {\frak X}^r(M)$,
arbitrarily close to $X,$ having
a closed trajectory passing through $p.$
\end{teo}

Related to this theorem (see [Gu2]), we have that:
For any  $E\in \BB$,
it can be constructed $Y\in {\frak X}^\infty(S),$
for some surface $S,$  
having a nontrivial recurrent point 
$p_0$ such that item (a1) is satisfied for some $h : [p_0,q_0) \to [0,1),$
and $P_Y :[p_0,q_0)\to [p_0,q_0).$ 
Here, $P_Y$ can be obtained to be injective or not.

\section{Proof of the results}

Suppose that $M$ is endowed with an orientation 
and with a
smooth riemannian metric $\,\, <\,\, , \,\,>. \, \,$ Given a
\xcrv we define $X^\perp \in {\frak X}^r(M)$ by the following 
conditions:
(a) \, \,  $<X,X> = <X^\perp, X^\perp>$; and \, \,
(b) \, \, when $p\in M$ is regular point of $X$, the 
ordered pair $(X(p), X^\perp(p))$ is an orthogonal positive
basis of $T_p(M)$ (according to the given orientation of $M$).
let  $\Sigma$ be an arc of trajectory of
$X^\perp.$  
A {\it \sfb \esp (for $X$)} is
a compact subset $F\subset M$ whose
interior is a flow box of $X$ and whose
boundary $\partial F$ is 
a  graph, homeomorphic to the
figure ``8'', which is the union of arcs of trajectory
$[\ov{c},\ov{a}]_X$ and $[\ov{a},\ov{c}]$  
(connecting $\ov{a}$ and $\ov{c}$) of  
$X$ and $X^\perp,$ respectively.
We shall refer to $[\ov{a},\ov{c}]$ (resp. $[\ov{c},\ov{a}]_X$)
as the {\it orthogonal \/} (resp. {\it tangent) edge \/} 
of either $F$ or $\partial F.$
See Figs. 1.a and 1.b.

\vspace{0.1cm}
% CHANGED BY JPZ \input{psfig}
\begin{figure}
\includegraphics[width=12cm,height=4.5cm]{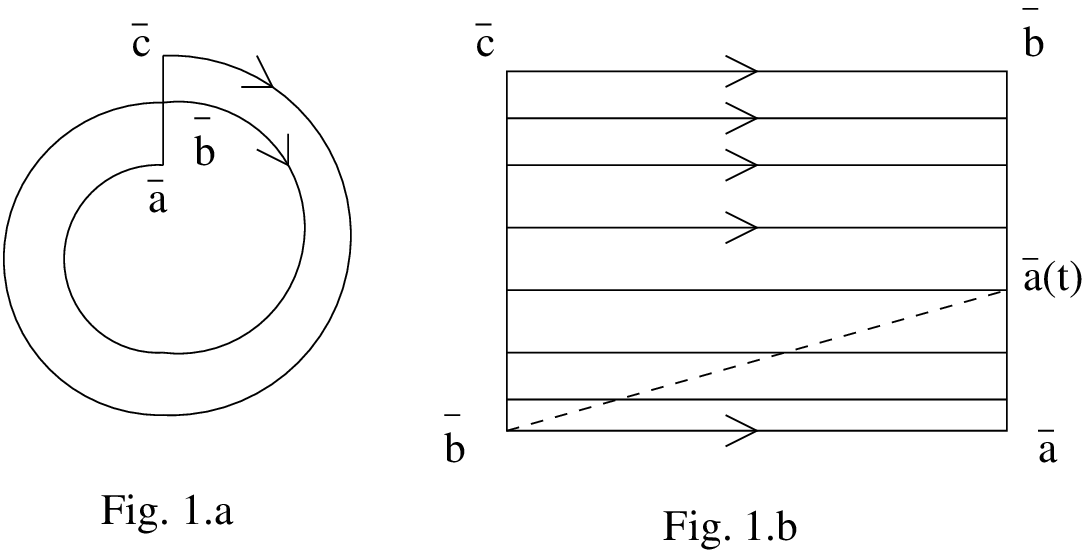} % CHANGED BY JPZ 
% CHANGED BY JPZ \centerline{\psfig{figure=1.eps,height=4.5cm,width=12cm}}
\end{figure}

%\vskip 4.5cm
%\centerline{Fig. 1.a \hskip 6cm Fig. 1.b}
 
%\bigskip

Let \xcr and let $p\in M$ be a 
nontrivial recurrent point of $X.$ 
We say that $X$ is {\it \tc \esp } at $p$ 
(i.e. twist-closable at $p$) if
there exists a half-open interval  $\Sigma=[p,q)$ 
tangent to $X^\perp,$
such that, for any neighborhood
$V$ of $p,$ there exists a \sfb \esp  for $X$ having 
its orthogonal edge contained in $\Sigma\cap V.$ 

\begin{prop} 
Let \xcr and let $p\in M$ be a 
nontrivial recurrent point of $X.$ 
Suppose that $X$ is \tc \esp at $p.$ 
Then there are 
sequences $t_n \to 0,$ of real numbers, and $p_n \to p,$
of points of $M,$  such that
$X + t_n X^\perp$ has a closed trajectory through $p_n$
\end{prop}

\no {\bf Proof:}
As $X$ is \tc \esp at $p,$ there exists a half-open interval $\Sigma=[p,q)$  
tangent to
$X^\perp,$ such that, 
given neighborhoods $\V$ of $X$ and $V$ of $p,$
we may choose a \sfb \esp $F\subset M$ (for $X$) and 
$\sigma > 0$ such that if
$[\ov{c},\ov{a}]_X$  and $[\ov{a},\ov{c}]$) are the
tangent and orthogonal edges, respectively, of $\partial F,$
and  
$\ov{b}$ is the vertex of $\partial F,$ then:
\begin{enumerate}
\item[{(b1)}] $[\ov{a},\ov{c}]\subset V$ 
and the flow of $X$ enters into $F$ through
the closed subinterval $[\ov{b},\ov{c}]$ of $\Sigma$;
moreover, for all $t\in [-\sigma, \sigma],$ 
$X(t):= X+ t \, X^\perp \in \V;$
\item[{(b2)}] both $X(\sigma)$ and $X(-\sigma)$ have
an arc of trajectory
contained in $F,$ which is a global cross section for 
$X\vert_F.$
\end{enumerate}

We shall continue considering only the case in which the
flow of $X^\perp$ goes from $\ov{a}$ to $\ov{c}.$
Let $\Gamma$ be the set of real numbers $s\in [0, \sigma]$ such that
when $t \in [0, s]$ there is an arc of trajectory 
$[\ov{b},\ov{a}(t)]_{X(t)}$ of $X(t),$ 
joining $\ov{b}$ with $\ov{a}(t)\in [\ov{a},\ov{b}]$, 
contained in $F,$ with $\ov{a}(0)= \ov{a},$ and 
such that $\ov{a}(t)$ depends 
continuously on $t.$ When $t\in \Gamma,$ 
these conditions determine 
$\ov{a}(t)$ and also that $[\ov{b},\ov{a}(t)]_{X(t)}$ is transversal to $X.$
Therefore, by (b2), $\Gamma = [0,\sigma_1]$ is  a
closed interval,
$\ov{a}(\sigma_1) = b$ and 
$[\ov{b}, \ov{a}(\sigma_1)]_{X(\sigma_1)}$
is a closed trajectory of $X{(\sigma_1)}.$
 See Fig. 1.b \qed

Under the assumptions and conclusions of this proposition,
there exists a sequence $F_n : M \to  M$ of $C^r-$diffeomorphisms, 
taking $p_n$ to $p$. We may assume that 
$F_n$ converges to the identity diffeomorphism in the $C^{r+1}-$topology.
Therefore, the sequence of
vector fields $(F_n)_* (X + t_n X^\perp) \to X$ in the $C^r-$topology
and each $(F_n)_* (X + t_n X^\perp)$ has a closed trajectory
passing through $p.$ This proves the following
\begin{teo}
Let \xcrp Let $p\in M$ be a nontrivial recurrent point of $X.$ 
Suppose that that $X$ is \tc \esp at $p.$ 
Then there exists $Y\in {\frak X}^r(M)$ arbitrarily close to $X$
having a closed trajectory through $p.$ 
\end{teo}

\medskip

\no {\bf Proof of Theorem 2.1: } 
We shall prove that: For all $m\ge 2,$ \,\, 
$\Delta_m \times {\frak S}_m \setminus \BB$
is a measure zero set.
It was proved by W. A. Veech [Vee] and H. Masur [Mas] 
that the Rauzy operator
$
{\cal R} : {\cal M} \to {\cal M},
$
defined in a full measure subset ${\cal M}$ of
$\Delta_m \times {\frak S}_m,$ is ergodic and 
has the following property:
\begin{enumerate}
\item[{(c)}] 
Given $E \in {\cal M},$  there exists a sequence $\{[0,a_n)\}$ of
subintervals  of $[0,1)$ such that $a_n \to 0$ and, if
$\tilde{E_n} : [0, a_n) \to [0, a_n)$ denotes the
\iet \esp induced by $E$, then, up to re-scaling, 
${\cal R}^n(E)$ coincides with $\tilde{E_n};$ more precisely,
${\cal R}^n(E) (z) = (1/a_n) \tilde{E}_n (a_n z)$, for
all $z\in [0,1)$.
\end{enumerate}

Given $k\ge 1$, let $A_k$ be the set of 
$E\in \Delta_m \times {\frak S}_m$
such that for some $a \in (16^{-k} - 32^{-k}, 
16^{-k} + 32^{-k}),$ $E(x) = a + x$, for all $x\in [0,1/2].$
We observe that $A_k$ is open and so it has positive
measure. 
Let $\tilde{\BB_k}$ be the set of $E\in {\cal M}$ such that
the positive ${\cal R}-$orbit of $E$ visits $A_k$ 
infinitely many often. As $A_k$ has positive measure and
${\cal R}$ is ergodic, the complement of  $\tilde{\BB_k}$ has measure
zero. Therefore, the complement of 
$
\tilde{\BB} = \cap_{k\ge 2} \tilde{\BB_k}
$
has measure zero.
Observe that if and \iet \esp $E\in A_k,$  then $E$ has 
more than $k,$ pairwise disjoint, \vtes. 
Therefore, as
${\cal R}$ satisfy  (c)  right above
and since the positive ${\cal R}-$orbit of
any given $E\in \tilde{\BB}$  visits every $A_k$ infinitely
many often, we obtain that $\tilde{\BB} \subset \BB.$ this 
proves the theorem. \qed  

\medskip

\no {\bf Proof of Theorem 2.3:}
This theorem is stated as follows:
Let $p\in M$ be a \bkt \esp nontrivial 
recurrent point of \xcrp
Suppose that $X$ has $K\ge 0$ singularities.
Then there exists a $Y\in {\frak X}^r(M)$
arbitrarily close to $X,$  having
a closed trajectory passing through $p.$

By theorem 3.2, 
it is enough to prove that $X$ is \tc \esp at $p.$
Let $\Sigma = [p,q),\,$ 
$T:[p,q)\to [p,q),\,$ $E:[0,1)\to [0,1),\,$
$h:[p,q)\to [0,1)\,$ be as in Theorem 2.2. 
As $E\in \BB_K,$  
given a neighborhood $V$ of $p,$ there exist
$b\in (0,1)$ and an \iet \esp
$E_V : [0, b) \to [0,b),$ such that:
\begin{enumerate}
\item[{(e)}] $E_V$ has at least $\chi+K+3$ 
pairwise disjoint \vte s contained in $[0,b);$ moreover,
the interval $\Sigma_V = h^{-1}([0,b))$ is contained in
$V.$
\end{enumerate}
Let $T_V : \Sigma_V \to \Sigma_V$ be the map induced by $T.$
As $X$ has $K$ singularities, (e) and
Theorem 2.2 imply that $E_V$ has a
\vte \esp $[a,E_V(a)]\subset [0,b)\, $  such that,
for some $\ov{a}\in \mbox{DOM}((T_{\Sigma_V})),\, \, $
$[\ov{a}, T_V(\ov{a})] = h^{-1}([a,E(a)]) 
\subset \mbox{DOM}(T\vert_{\Sigma_V}).\, $
Therefore,  
there exists a \sfb \esp 
bounded by 
$[\ov{a}, T_V^2(\ov{a})] \cup [\ov{a}, T_V^2(\ov{a})]_X.$
As $V$ is arbitrary,
this proves that $X$ is \tc \esp at $p.$ \qed

\end{document}